\documentclass[A4paper,11pt]{article}
\usepackage{amsfonts}
\usepackage{amsmath}
\usepackage{amssymb}
\usepackage{graphicx}
\begin{document}
\title{Properties of Khovanov homology for positive braid knots  
\author{Marko Sto\v si\'c 
\thanks{The author is supported by {\it Funda\c c\~ao de Ci\^encia e
Tecnologia}/(FCT), grant no. SFRH/BD/6783/2001}\\
Departamento de Matem\'atica  and \\
CEMAT -  Centro de Matem\'atica e Aplica\c c\~oes\\
Instituto Superior T\'ecnico\\
Av. Rovisco Pais 1\\
1049-001 Lisbon\\ 
Portugal\\
e-mail: mstosic@math.ist.utl.pt
}
}

\date{}

\newtheorem{theorem}{Theorem}
\newtheorem{acknowledgment}[theorem]{Acknowledgment}
\newtheorem{algorithm}[theorem]{Algorithm}
\newtheorem{axiom}[theorem]{Axiom}
\newtheorem{case}[theorem]{Case}
\newtheorem{claim}[theorem]{Claim}
\newtheorem{conclusion}[theorem]{Conclusion}
\newtheorem{condition}[theorem]{Condition}
\newtheorem{conjecture}[theorem]{Conjecture}
\newtheorem{corollary}[theorem]{Corollary}
\newtheorem{criterion}[theorem]{Criterion}
\newtheorem{definition}{Definition}
\newtheorem{example}[theorem]{Example}
\newtheorem{exercise}[theorem]{Exercise}
\newtheorem{lemma}{\indent Lemma}
\newtheorem{notation}[theorem]{Notation}
\newtheorem{problem}[theorem]{Problem}
\newtheorem{proposition}{Proposition}
\newtheorem{remark}[theorem]{Remark}
\newtheorem{solution}[theorem]{Solution}
\newtheorem{summary}[theorem]{Summary}
\newcommand{\ud}{\mathrm{d}}

\def\gcd{\mathop{\rm gcd}}
\def\Ker{\mathop{\rm Ker}}
\def\max{\mathop{\rm max}}
\def\map{\mathop{\rm map}}
\def\lcm{\mathop{\rm lcm}}
\def\kraj{\hfill\rule{6pt}{6pt}}
\def\diag{\mathop{\rm diag}}
\def\span{\mathop{\rm span}}
\def\deg{\mathop{\rm deg}}
\def\rank{\mathop{\rm rank}}
\def\sgn{\mathop{\rm sgn}}
\def\kvn{1+q+\ldots+q^{n-1}}
\def\F{\mathbb{F}}
\def\R{\mathbb{R}}
\def\C{\mathcal{C}}
\def\K{\mathbb{K}}
\def\Z{\mathbb{Z}}
\def\Q{\mathbb{Q}}
\def\H{\mathcal{H}}
\arraycolsep 6pt

\maketitle

\begin{abstract} 
In this paper we solve one open problem from \cite{pat} and give some generalizations. Namely, we prove that the first homology group of positive braid knot is trivial. Also, we show that the same is true for the Khovanov-Rozansky homology \cite{kovroz} ($sl(n)$ link homology) for any positive integer $n$. 
\end{abstract}

\section{Introduction}

\indent 
In recent years there has been a lot of interest in the 
 ``categorification" of link invariants, initiated by Khovanov in \cite{kov}. For each link $L$ in $S^3$ he defined a graded chain complex, with grading 
preserving differentials, whose graded Euler characteristic is equal to
the Jones polynomial of the link $L$ (\cite{jones},\cite{kauf}).  This is done by starting from the state sum 
expression for the Jones polynomial (which is written as an alternating sum),
then constructing for each term a module whose graded dimension is 
equal to the value of that term, and finally, constructing the
differentials 
as appropriate grading preserving maps, so that 
the complex obtained is a link invariant. There
is also similar construction \cite{kovroz},  for the categorification of the $n$-specializations of HOMFLYPT polynomial (\cite{homfly}, \cite{pt}, \cite{moy}) as well as the categorifications of various 
link and graph polynomial invariants (\cite{kov3}, \cite{kol}, \cite{bn1}, \cite{gr}).\\
\indent Although the theory is rather new, it already has strong applications. For instance, the proof of Milnor conjecture by Rasmussen in \cite{ras} and the existence of the exotic differential structure on
 $\R^4$ (\cite{ras2},\cite{xy} page 522), which were previously accesible only by huge machinery of gauge theory.\\
\indent     The advantage of Khovanov homology theory is that its definition is combinatorial and since there is a straightforward algorithm for computing it, it is (theoretically) highly calculable. Nowadays there are several computer programs \cite{bnprog}, \cite{shum} that  can calculate effectively Khovanov homology of links with up to 50 crossings.\\
\indent Based on the calculations there are many conjectures about the properties of link homology see e.g. \cite{bn}, \cite{pat}, \cite{dgr}. Some of the properties have been proved till now (see \cite{lee}, \cite{lee2}), but many of them are still open. \\
\indent In this paper we prove the conjecture from \cite{pat} that the first homology group of the positive braid knot is trivial. In the proof we use the basic ingredients of the construction of homology (cubic complex, independence of the planar projection chosen), and so major part can be directly applied to other link homology theories. Especially, in the case of $sl(n)$-link homology \cite{kovroz}, whose particular definitions of chain groups and differentials make it practically incalculable, we modify slightly our proof to show that the first $sl(n)$ homology group  of positive braid knot is trivial, as well.\\
\indent Since torus knots are positive braid knots, we automatically obtain the same properties for them. Even more, in the sequel papers \cite{mojtorus}, \cite{mojnov}, we obtain further properties of $sl(2)$   
 and $sl(n)$-link homology for torus knots.\\
\indent The organization of the paper is the following: in Section \ref{intro} we briefly recall the definition of the Khovanov homology. In Section \ref{glavni} we give the main result. Finally, in Section \ref{sln} we adapt the proof of the main result for $sl(n)$ case.\\

\textbf{Acknowledgements:}
The author would like to thank M. Khovanov and J. Rasmussen for many helpful discussions.

\section{Notation}\label{intro}

We recall briefly the definition of Khovanov homology for links. For more details see \cite{bn},\cite{kov}.  \\
\indent First of all, take a link $K$, its planar projection $D$, and take an ordering of the crossings of $D$. For each crossing $c$ of $D$, we define 0-resolution and 1-resolution, according to the following picture:
 
{{
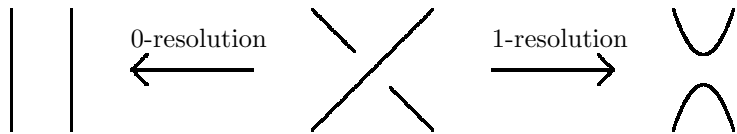
\begin{center}
\setlength{\unitlength}{4mm}
\begin{picture}(30,6)
\linethickness{0.8pt}
\qbezier(1.00,1.00)(1.00,3.00)(1.00,5.00)
\qbezier(3.00,1.00)(3.00,3.00)(3.00,5.00)

\qbezier(5.00,3.00)(7.00,3.00)(9.00,3.00)
\qbezier(5.00,3.00)(5.25,3.25)(5.50,3.50)
\qbezier(5.00,3.00)(5.25,2.75)(5.50,2.50)

\qbezier(15.00,5.00)(13.00,3.00)(11.00,1.00)
\qbezier(15.00,1.00)(14.30,1.70)(13.60,2.40)
\qbezier(11.00,5.00)(11.70,4.30)(12.40,3.60)

\qbezier(17.00,3.00)(19.00,3.00)(21.00,3.00)
\qbezier(21.00,3.00)(20.75,3.25)(20.50,3.50)
\qbezier(21.00,3.00)(20.75,2.75)(20.50,2.50)

\qbezier(23.00,5.00)(24.00,2.00)(25.00,5.00)
\qbezier(23.00,1.00)(24.00,4.00)(25.00,1.00)

\put(5.00,3.75){{\small\textrm{0-resolution}}}
\put(17.00,3.75){{\small\textrm{1-resolution}}}

\end{picture}
\end{center}
}}

Denote by $m$ the number of crossings of  $D$. Then there is bijective correspondence between the total resolutions of $D$ and the set $\{0,1\}^m$. Namely, to every $m$-tuple $\epsilon=(\epsilon_1,\ldots,\epsilon_m)\in \{0,1\}^m$ we associate the resolution $D_{\epsilon}$ where we resolved the $i$-th crossing in a $\epsilon_i$-resolution.\\
\indent Every resolution $D_{\epsilon}$ is a collection of disjoint circles. To each circle we associate graded $\Z$-module $V$, which is freely generated by two basis vectors $1$ and $X$, with $\deg 1 = 1$ and $\deg X=-1$. To $D_\epsilon$ we associate the module $M_{\epsilon}$, which is the tensor product of $V$'s over all circles in the resolution. We group all the resolution $D_\epsilon$ with fixed $|\epsilon|$ (sum of elements of $\epsilon$). We draw all resolutions as (skewed) $m$-dimensional cube such that in $i$-th column are the resolutions $D_{\epsilon}$ with $|\epsilon|=i$. We define the $i$-th chain group $C^i$ by:
$$C^i(D)=\oplus_{|\epsilon|=i}M_{\epsilon}\{i\}.$$
\indent Here, by $\{i\}$, we have denoted shift in grading of $M_{\epsilon}$ (for more details see e.g. \cite{bn}).\\
\indent The differential $d^i:C^i(G)\to C^{i+1}(G)$ is defined as (signed) sum of  ``per-edge" differentials. Namely, the only nonzero maps are from $D_{\epsilon}$ to $D_{\epsilon'}$, where $\epsilon=(\epsilon_1,\ldots,\epsilon_m)$, $\epsilon_i \in \{0,1\}$, if and only if $\epsilon'$ has all entries same as $\epsilon$ except one $\epsilon_j$, for some $j\in \{0,1\}$, which is changed from 0 to 1. We denote these differentials by $d_{\nu}$, where $\nu$ is $m$-tuple which consists of the label $\ast$ at the position $j$ and of $m-1$ 0's and 1's (the same as the remaining entries of $\epsilon$). 
Note that in these cases, either two circles of $D_{\epsilon}$ merge into one circle of $D_{\epsilon'}$ or one circle of $D_{\epsilon}$ splits into two circles of $D_{\epsilon'}$, and all other circles remain the same. \\
\indent In the first case, we define the map $d_{\nu}$ as the identity on the tensor factors ($V$) that correspond to the unchanged circles, and on the remaining factors we define as the (graded preserving) multiplication map $m:V\otimes V \to V\{1\}$, which is given on basis vectors by:
$$m(1\otimes 1)=1,\quad m(1\otimes X)=m(X\otimes 1)=X,\quad m(X\otimes X)=0.$$
\indent In the second case, we define the map $d_{\nu}$ as the identity on the tensor factors ($V$) that correspond to the unchanged circles,
and on the remaining factors we define as the (graded preserving) comultiplication map $\Delta:V\to V \otimes V\{1\}$, which is given on basis vectors by:
$$\Delta(1)=1\otimes X+X\otimes 1,\quad \Delta(X)=X\otimes X.$$
\indent Finally, to obtain the differential $d^i$ of the chain complex $C(D)$, we sum all contributions $d_{\nu}$ with $|\nu|=i$, multiplied by the sign $(-1)^{f(\nu)}$, where $f(\nu)$ is equal to the number of 1's ordered before $\ast$ in $\nu$. This makes every square of our cubic complex anticommutative, 
 we obtain the genuine differential (i.e. $(d^i)^2=0$).\\
\indent The homology groups of the obtained complex $(C(D),d)$ 
we denote by $H^i(D)$ 
and call \textit{unnormalized homology groups of $D$}. In order to obtain link invariants (i.e.
independence of the chosen projection), we have to shift the chain complex (and hence the homology groups) by:
\begin{equation}
\C(D)=C(D)[-n_-]\{n_+ -2n_- \},
\label{eq1}\end{equation}
where $n_+$ and $n_-$ are number of positive and negative crossings, respectively, of the diagram $D$. 
{{\begin{center}
\setlength{\unitlength}{4mm}
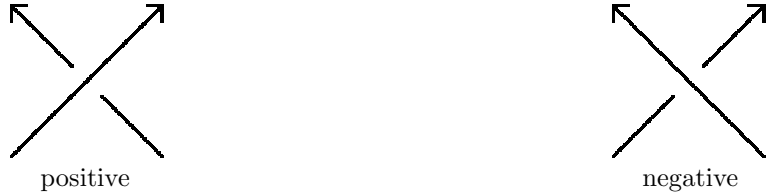
\begin{picture}(25,6)
\linethickness{0.8pt}
%\line(1.00,2.00,1.00,3.00)
\qbezier(5.00,5.00)(2.50,2.50)(0.00,0.00)
\qbezier(5.00,0.00)(4.00,1.00)(3.00,2.00)
\qbezier(0.00,5.00)(1.00,4.00)(2.00,3.00)
\qbezier(0.00,5.00)(0.25, 5.00)(0.50,5.00)
\qbezier(0.00,5.00)(0.00,4.75)(0.00,4.50)
\qbezier(5.00,5.00)(4.75, 5.00)(4.50,5.00)
\qbezier(5.00,5.00)(5.00,4.75)(5.00,4.50)
\qbezier(22.00,2.00)(21.00,1.00)(20.00,0.00)
\qbezier(25.00,5.00)(24.00,4.00)(23.00,3.00)
\qbezier(20.00,5.00)(22.50,2.50)(25.00,0.00)
\qbezier(20.00,5.00)(20.25,5.00)(20.50,5.00)
\qbezier(20.00,5.00)(20.00,4.75)(20.00,4.50)
\qbezier(25.00,5.00)(24.75,5.00)(24.50,5.00)
\qbezier(25.00,5.00)(25.00,4.75)(25.00,4.50)
\put(1.00,-1.00){{\small\textrm{positive}}}
\put(21.00,-1.00){{\small\textrm{negative}}}
%\put(10.00,-2.00){{\small\textrm{figure 1}}}
\end{picture}
\end{center}
}}
\bigskip
In the formula (\ref{eq1}), we have denoted by $[-n_-]$, the shift in homology degrees (again, for more details see \cite{bn}).\\
\indent The homology groups of the complex $\C(D)$ we denote by $\H^i(D)$. Hence, we have $\H^{i,j}(D)=H^{i+n_-,j-n_++2n_-}(D)$. 

\begin{theorem}(\cite{kov},\cite{bn})
The homology groups $\H(D)$ are independent of the choice of the planar projection $D$. Furthermore, the graded Euler characteristic of the complex $\C(D)$ is equal to Jones polynomial of the link $K$.
\end{theorem}
\indent Hence, we can write $\H(K)$, and we call $\H^i(K)$ \textit{homology groups of link $K$}.\\

\indent Also if the diagram $D$ has only positive crossings then we do not have shift in the homology degrees, and so we have that, for example, $\H^1(K)$ is trivial if and only if $H^{1}(D)$ is trivial. Also, in general case if we have a positive knot $K$ (the knot that has a planar  projection with only positive crossings) then we have that $\H^i(K)$ is trivial for all $i<0$. Further, if $D$ is planar projection of positive knot $K$ with $n_-$ negative crossings then $H^i(D)$ is trivial for $i<n_-$. \\

\section{Positive braid knots}\label{glavni}
  
The positive braid knots are the knots (or links) that are the closures of positive braids. Let $K$ be arbitrary positive braid knot and let $D$ be its planar projection which is the closure of a positive braid. Denote the number of strands of that braid by $p$. We say that  the crossing $c$ of $D$ is of the type $\sigma_i$, $i<p$, if it corresponds to the generator $\sigma_i$ in the braid word of which $D$ is the closure.

{{
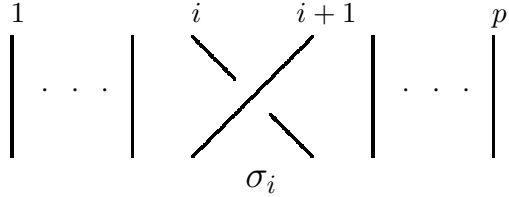
\begin{center}
\setlength{\unitlength}{4mm}
\begin{picture}(30,6) 
\linethickness{0.8pt}
%\line(1.00,2.00,1.00,3.00)
\qbezier(5.00,1.00)(5.00,3.00)(5.00,5.00)
\qbezier(9.00,1.00)(9.00,3.00)(9.00,5.00)

\qbezier(15.00,5.00)(13.00,3.00)(11.00,1.00)
\qbezier(15.00,1.00)(14.30,1.70)(13.60,2.40)
\qbezier(11.00,5.00)(11.70,4.30)(12.40,3.60)

\qbezier(17.00,1.00)(17.00,3.00)(17.00,5.00)
\qbezier(21.00,1.00)(21.00,3.00)(21.00,5.00)

\put(6.00,3.00){$\cdot$}
\put(7.00,3.00){$\cdot$}
\put(8.00,3.00){$\cdot$}
\put(18.00,3.00){$\cdot$}
\put(19.00,3.00){$\cdot$}
\put(20.00,3.00){$\cdot$}

\put(5.00,5.50){1}
\put(11.00,5.50){$i$}
\put(14.50,5.50){$i+1$}
\put(21.00,5.50){$p$}
\put(12.80,0.00){{\Large $\sigma_i$}}

\end{picture}
\end{center}
}}

 Denote the number of crossings of the type $\sigma_i$ by $l_i$, $i=1,\ldots,p-1$ and order them from top to bottom. Then each crossing $c$ of $D$ we can write as the pair $(i,\alpha)$ (we will also write $(i\alpha)$ if there is no possibility of confusion), $i=1,\ldots,p-1$ and $\alpha=1,\ldots,l_i$, if $c$ is of the type $\sigma_i$ and it is ordered as $\alpha$-th among the crossings of the type $\sigma_i$. Finally, we order the crossings of $D$ by the following ordering: $c=(i\alpha)<d=(j\beta)$ if and only if $i<j$, or $i=j$ and $\alpha<\beta$.\\
\indent Since the positive braid knot is positive knot, we have that $\H^i(K)=0$ for $i<0$. Also we know its zeroth homology group (see e.g. \cite{pat}) is two-dimensional (without torsion) and that the $q$-gradings of the two generators are $1-p+n(D)\pm 1$, where $n(D)$ is the number of crossings of $D$. In the following theorem we prove that the first homology group of the positive braid knot is trivial.
\begin{theorem}
If $K$ is positive braid knot, than $\H^1(K)=0$.
\end{theorem}
\textbf{Proof:}\\
\indent First of all, if $D$  is the regular diagram of $K$, which is the closure of a positive braid, then we have that $\H^{i,j}(K)=H^{i,j-n(D)}(D)$, where $n(D)$ is the number of crossings of $D$, and so we have that $\H^1(K)=0$ if and only if $H^1(D)=0$.  So we are going to show that the (unnormalized) first homology group of $D$ is trivial. \\
\indent In order to prove that $H^1(D)=0$, we will use the definition of Khovanov homology, i.e. we will prove that for any element $t'$ from chain group $C^1(D)$ such that $d^1(t')=0$ there exists an element $t\in C^0(D)$ such that $t'=d^0(t)$. For this we first need to understand the chain groups $C^0(D)$, $C^1(D)$ and $C^2(D)$ and the differentials  $d^0$ and $d^1$.\\
\indent $C^0(D)$ ``comes" from all resolutions $K_s$ with $|s|=0$. However, since $D$ is the closure of the positive braid we have only one such resolution $s_0$ (all crossings are resolved into 0-resolutions) and it is  an unlink  that consists  of $p$ unknots (i.e. the closure of the trivial braid with $p$ strands). Hence, we have that $C^0(D)=V^{\otimes p}$ where we assigned the  $i$-th  copy of $V$ (denoted by $V^i$) to the circle  which is the closure of the $i$-th strand of $D_{s_0}$. \\
\indent Now, we pass to $C^1(D)$. It ``comes"  from all  resolutions $D_s$ with $|s|=1$, i.e. all the resolutions were we resolve all except one crossing of $D$ in a 0-resolutions  
 and the remaining one in 1-resolution. 
In such way, if the crossing $c$ that is resolved into 1-resolution is of the type $\sigma_i$ (i.e. if $c=(i\alpha)$, for some $\alpha=1,\ldots,l_i$), then the corresponding resolution, $D_c$, is the closure of the plat diagram $E_i$ (see picture).

{{
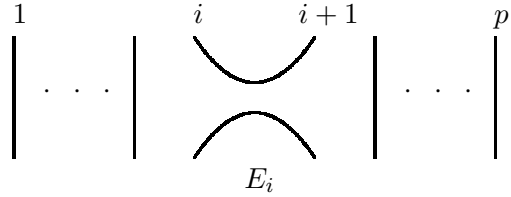
\begin{center}
\setlength{\unitlength}{4mm}
\begin{picture}(30,6)
\linethickness{0.8pt}
%\line(1.00,2.00,1.00,3.00)
\qbezier(5.00,1.00)(5.00,3.00)(5.00,5.00)
\qbezier(9.00,1.00)(9.00,3.00)(9.00,5.00)

\qbezier(15.00,5.00)(13.00,2.00)(11.00,5.00)
\qbezier(11.00,1.00)(13.00,4.00)(15.00,1.00)

\qbezier(17.00,1.00)(17.00,3.00)(17.00,5.00)
\qbezier(21.00,1.00)(21.00,3.00)(21.00,5.00)

\put(6.00,3.00){$\cdot$}
\put(7.00,3.00){$\cdot$}
\put(8.00,3.00){$\cdot$}
\put(18.00,3.00){$\cdot$}
\put(19.00,3.00){$\cdot$}
\put(20.00,3.00){$\cdot$}

\put(5.00,5.50){1}
\put(11.00,5.50){$i$}
\put(14.50,5.50){$i+1$}
\put(21.00,5.50){$p$}
\put(12.70,0.00){$E_i$}

\end{picture}
\end{center}
}}

To that resolution we assign the vector space $V_c=V^{\otimes(p-1)}$, where we have assigned the first $i-1$ and the last $p-1-i$ copies of $V$ to the circles that are the closures of the first $i-1$ and the last $p-i$ strands of the resolution $D_c$, respectively, and the $i$-th copy of $V$ corresponds to the remaining circle (closure of $E_i$). So, we have that $C^1(D)=\bigoplus_{c\in c(D)}{V_c \{1\}}$. Further on, we denote the $k$-th copy of $V$ in $V_c$ by $V^k_c$.\\ 

\indent The differential $d^0: C^0(D)\to C^1(D)$ is given by the maps $f_{c}: V^{\otimes p}\to V_{c}$ where if $c$ is of the type $\sigma_i$ than $f_c$ acts as the identity on the first $i-1$ and the last $p-i-1$ copies of $V$ and as the multiplication $m$ on the remaining two copies of $V$. In other words, $f_c$ maps the copies $V^j$ as the identity onto $V^j_c$, for $j<i$, maps the copies $V^{j+1}$ as the identity onto $V^{j}_c$, for $i<j<p$, and on the remaining two factors act as the  multiplication $m:V^i\otimes V^{i+1} \to V^i$.\\

\indent Further, $C^2(D)$ comes from the resolutions where exactly two of the crossings are resolved in 1-resolutions and the remaining ones are resolved in 0-rezolutions. Denote the two  crossings that are resolved in  1-resolution  by $c$ and $d$ (where $c$ is ordered before $d$), and let $c$ be of $\sigma_i$ type and $d$ of $\sigma_j$ type. Then we have that $i \le j$.\\
\indent If $i=j$ then the corresponding resolution $D_{c,d}$ is a closure of a plat diagram $E_i^2$, has $p$ circles and hence the corresponding summand  $V_{c,d}$ of a $C^2(D)$ is isomorphic to $V^{\otimes p}$. Here we have assigned the first $i-1$ and the last $p-i-1$ copies of $V$ to the circles that are the closures of the first $i-1$ and the last $p-i-1$ strands of the resolution $D_{c,d}$, respectively, while  $i$-th and $(i+1)$-th copy of $V$ correspond to the remaining two circle that are formed of $i$-th and $(i+1)$-th strand ($i$-th copy of $V$ to the outer, and $(i+1)$-th copy of $V$ to the inner circle). \\
\indent If $i<j$ then the corresponding resolution is a closure of a plat diagram $E_i E_j$ (or $E_j E_i$), has $p-2$ circles and hence the corresponding summand $V_{c,d}$ of a $C^2(D)$ is isomorphic to $V^{\otimes (p-2)}$. \\
\indent If $i+1<j$, we assign the first $i-1$ copies of $V$ to the closures of the first $i-1$ strands. The $i$-th copy of $V$ is assigned to the circle that is obtained by joining the $i$-th and $(i+1)$-th strand (the 1-resolution of $c$). The following $j-i-2$ copies of $V$ are assigned to the closures of the strands from $(i+2)$-th to $(j-1)$-th of $D_{c,d}$, respectively.
The $(j-1)$-th copy of $V$ is assigned to the circle that is obtained by joining the $j$-th and $(j+1)$-th strand (the 1-resolution of $d$). 
The remaining $p-j-1$ copies of $V$ are assigned to the closures of the last $p-j-1$ strands of $D_{c,d}$. \\
\indent If $i+1=j$, then we assign the first $i-1$ copies of $V$ to the closures of the first $i-1$ strands. The $i$-th copy of $V$ is assigned to the circle that is obtained by joining the $i$-th, $(i+1)$-th and $(i+2)$-th strand (the 1-resolutions of $c$ and $d$). The remaining $p-i-2$ copies of $V$ are assigned to the closures of the strands from $(i+3)$-th to $p$-th of $D_{c,d}$, respectively.\\
\indent In all previous cases, we denote the $k$-th copy of $V$ in $V_{c,d}$ by $V^k_{c,d}$.\\
 
\indent Finally, the second chain group is $$C^2(D)=\bigoplus_{c,d\in c(D),\,c<d} V_{c,d}\{2\}.$$\\

\indent Now, we can describe the differential $d^1: C^1(D)\to C^2(D)$. 
It is given as a sum of the maps of the form $f_{ecd}: V_e\to V_{c,d}$ (with $c,d,e\in c(D),\,c<d$), where $f_{ecd}$ is zero unless $e=c$ or $e=d$. Let  $c=(i\alpha)$ and $d=(j\beta)$, for some $1\le i \le j \le n-1$, $\alpha=1,\ldots,l_i$ and $\beta=1,\ldots,l_j$. The maps $f_{ccd}: V_c\to V_{c,d}$ are as follows:  \\
\indent If $i=j$ then $f_{ccd}$ is given by the identity maps: $id:V^l_c\to V^l_{c,d}$, for $l<i$, and $id:V^{l}_c\to V^{l+1}_{c,d}$, for $i<l<p$, and by comultiplication $\Delta:V^i_c\to V^i_{c,d}\otimes V^{i+1}_{c,d}$.\\
\indent If $i<j$, then $f_{ccd}$ is given by the identity maps
$id:V^l_c\to V^l_{c,d}$, for $l<j-1$, and $id:V^{l+1}_c\to V^{l}_{c,d}$, for $j<l<p-1$, and by multiplication $m:V^{j-1}_c\otimes V^j_c \to V^{j-1}_{c,d}$.\\
\indent The other class of nonzero maps $f_{dcd}: V_d\to V_{c,d}$ is in the case $i=j$ given by $-f_{ccd}$ (since in this case $V_c=V_d$, and $c<d$), and in the case $i<j$ is given as $-g_{cd}$, where $g_{cd}$ is given by the identity maps:
$id:V^l_c\to V^l_{c,d}$, for $l<i$, and $id:V^{l+1}_c\to V^{l}_{c,d}$, for $i<l<p-1$, and by multiplication $m:V^{i}_c\otimes V^{i+1}_c \to V^{i}_{c,d}$.\\ 

\indent
Now, we can go back to the proof. Let $$t'=(t_{1,1},\ldots,t_{1,l_1},t_{2,1},\ldots,t_{2,l_2},\ldots,t_{p-1,1},\ldots,t_{p-1,l_{p-1}})\in C^1(D)$$ be such that $d^1(t')=0$. Here we have that  $t_{i,\alpha}\in V_{(i\alpha)}$ for $i=1,\ldots,p-1$, $\alpha=1,\ldots,l_i$. Our aim is to find an element $t\in C^0=V^{\otimes p}$ such that $d^0(t)=t'$. \\
\indent Since $d^1(t')=0$ we have that its projection, denoted by $d^1_{i\alpha\beta}$, to the space $V_{i\alpha,i\beta}$  is equal to zero, for every $i=1,\ldots,p-1$, and $\alpha,\beta=1,\ldots, l_i$,. 
However, this implies that $t_{i\alpha}=t_{i\beta}$ for every
$i=1,\ldots,p-1$, $\alpha,\beta=1,\ldots, l_i$, since only the maps from $V_{(i\alpha)}$ and $V_{(i\beta)}$ to $V_{(i\alpha),(i\beta)}$ are nonzero, and the only nonidentity part of the mappings is the comultiplication $\Delta$ on the same ($i$-th) copy of $V$ in both $V_{i\alpha}$ and $V_{i\beta}$. \\
\indent 
Hence, we have obtained that $t'\in\ker d^1$ if and only if     
$t_{i\alpha}=t_{i\beta}$ for every
$i=1,\ldots,p-1$, $\alpha,\beta=1,\ldots, l_i$ and $\bar{t}=(t_{1,1},t_{2,1},\ldots,t_{p-1,1})\in\ker\bar{d}^1$, where $\bar{d}^1$ is the restriction of $d^1$ to $W=V_{(1,1)}\oplus V_{(2,1)}\oplus\cdots V_{(p-1,1)}$.\\
\indent Further, note that for every $i=1,\ldots,p-1$, and $\alpha,\beta=1,\ldots, l_i$, the projection of $d^0(t)$, for any $t\in C^0(D)$ to $V_{(i\alpha)}$ and $V_{(i\beta)}$ is equal, i.e. the differential $d^0$ is completely determined by the map $\bar{d}^0$ which is the projection of $d^0$ on $W$. Finally, we have that 
if there exists $t\in C^0$ such 
$\bar{d}^0(t)=\bar{t}$, 
then $d^0(t)=t'$. Hence, to finish the proof, we are left with proving that for every $y\in \ker \bar{d}^1 \in W$, there exists $x \in C^0(D)$, such that $\bar{d}^0(x)=y$.\\
  
\indent 
Now observe the positive braid  knot $K'$ which has the regular diagram $D'$ which is the closure of the braid $\sigma_1\sigma_2\cdots\sigma_{p-1}$ (we omit the letters $\sigma_i$'s which are not contained in the braid word of which $D$ is the closure).  Its zeroth chain group $\bar{C}^0(D')$ is obviously equal to $C^0(D)=V^{\otimes p}$, the first chain group $\bar{C}^1(D')$ is equal to $W$ and its second chain group $\bar{C}^2(D)$ is equal to $\bigoplus_{1\le i<j\le p-1}{V_{(i1),(j1)}}$. Its zeroth differential is equal to the previously defined $\bar{d}^0$, while the first differential is equal to $\bar{d}^1$. Since $K'$ is isotopic to the unknot (or to the unlink consisting of the unknots), its first homology group is trivial and hence for every $y\in\ker \bar{d}^1$ there exists $x\in \bar{C}^0(D')=C^0(D)$ such that $\bar{d}^0(x)=y$. This concludes the proof. \kraj \\

\section{$sl(n)$ case}\label{sln}

In this section we will adapt our proof from the previous section to show that the 
first Khovanov-Rozansky homology group of a positive braid knot is trivial. We will not define Khovanov-Rozansky $sl(n)$-link
homology, for details see \cite{kovroz}. For our purposes it is enough to point out the similarities and the differences
between $sl(n)$-link homology and the standard Khovanov ($sl(2)$) homology from Section \ref{intro}.\\

\indent The basic principle of the construction, namely the cubic complex, is the same for $sl(n)$ case as it is in the
standard Khovanov homology. The difference is that we take the oriented knot and the two (oriented) resolutions of the
crossings as given in the following picture: \\

{{
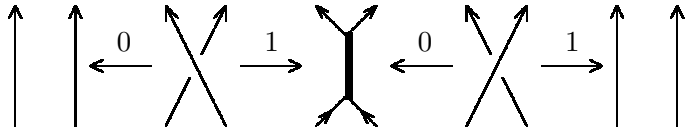
\begin{center}
\setlength{\unitlength}{4mm}
\begin{picture}(30,6) 
\linethickness{0.6pt}
%\line(1.00,2.00,1.00,3.00)
\qbezier(4.00,1.00)(4.00,3.00)(4.00,5.00)
\qbezier(3.80,4.60)(3.90,4.80)(4.00,5.00)
\qbezier(4.20,4.60)(4.10,4.80)(4.00,5.00)

\qbezier(2.00,1.00)(2.00,3.00)(2.00,5.00)
\qbezier(1.80,4.60)(1.90,4.80)(2.00,5.00)
\qbezier(2.20,4.60)(2.10,4.80)(2.00,5.00)

\qbezier(9.5,3)(10.5,3)(11.5,3)
\qbezier(11.5,3)(11.3,3.1)(11.1,3.2)
\qbezier(11.5,3)(11.3,2.9)(11.1,2.8)

\qbezier(4.5,3)(5.5,3)(6.5,3)
\qbezier(4.5,3)(4.7,3.1)(4.9,3.2)
\qbezier(4.5,3)(4.7,2.9)(4.9,2.8)

\qbezier(9,1)(8,3)(7,5)
\qbezier(7,1)(7.4,1.8)(7.8,2.6)
\qbezier(9,5)(8.6,4.2)(8.2,3.4)

\linethickness{2.5pt}
\qbezier(13.00,2.00)(13.00,3.00)(13.00,4.00)
\linethickness{0.6pt}
\qbezier(12.00,1.00)(12.50,1.50)(13.00,2.00)
\qbezier(14.00,1.00)(13.50,1.50)(13.00,2.00)
\qbezier(12.00,5.00)(12.50,4.50)(13.00,4.00)
\qbezier(14.00,5.00)(13.50,4.50)(13.00,4.00)

\qbezier(14.00,5.00)(13.80,4.90)(13.60,4.80)
\qbezier(14.00,5.00)(13.90,4.80)(13.80,4.60)

\qbezier(12.00,5.00)(12.20,4.90)(12.40,4.80)
\qbezier(12.00,5.00)(12.10,4.80)(12.20,4.60)

\qbezier(12.50,1.50)(12.30,1.40)(12.10,1.30)\qbezier(12.50,1.50)(12.40,1.30)(12.30,1.10)

\qbezier(13.50,1.50)(13.70,1.40)(13.90,1.30)\qbezier(13.50,1.50)(13.60,1.30)(13.70,1.10)

\qbezier(14.5,3)(15.5,3)(16.5,3)
\qbezier(14.5,3)(14.7,3.1)(14.9,3.2)
\qbezier(14.5,3)(14.7,2.9)(14.9,2.8)

\qbezier(19.5,3)(20.5,3)(21.5,3)
\qbezier(21.5,3)(21.3,3.1)(21.1,3.2)
\qbezier(21.5,3)(21.3,2.9)(21.1,2.8)

\qbezier(19,5)(18,3)(17,1)
\qbezier(17,5)(17.4,4.2)(17.8,3.4)
\qbezier(19,1)(18.6,1.8)(18.2,2.6)

\qbezier(19,5)(19,4.75)(19,4.5)
\qbezier(19,5)(18.8,4.8)(18.6,4.6)
\qbezier(17,5)(17,4.75)(17,4.5)
\qbezier(17,5)(17.2,4.8)(17.4,4.6)

\qbezier(9,5)(9,4.75)(9,4.5)
\qbezier(9,5)(8.8,4.8)(8.6,4.6)
\qbezier(7,5)(7,4.75)(7,4.5)
\qbezier(7,5)(7.2,4.8)(7.4,4.6)

%\qbezier(17.00,1.00)(17.00,3.00)(17.00,5.00)\qbezier(16.80,4.60)(16.90,4.80)(17.00,5.00)\qbezier(17.20,4.60)(17.10,4.80)(17.00,5.00)

\qbezier(22.00,1.00)(22.00,3.00)(22.00,5.00)
\qbezier(21.80,4.60)(21.90,4.80)(22.00,5.00)
\qbezier(22.20,4.60)(22.10,4.80)(22.00,5.00)

\qbezier(24.00,1.00)(24.00,3.00)(24.00,5.00)
\qbezier(23.80,4.60)(23.90,4.80)(24.00,5.00)
\qbezier(24.20,4.60)(24.10,4.80)(24.00,5.00)

\put(5.40,3.50){0}
\put(10.3,3.50){1}
\put(15.40,3.50){0}
\put(20.3,3.50){1}

\end{picture}
\end{center}
}}

\indent Hence, in this case the (complete) resolution of the planar projection $D$ consists of the three-valent graphs such
that at each vertex we have exactly one (unoriented) wide edge and two oriented thin edges, and such that at one end of each
wide edge two thin edges are incoming and at the other end are outgoing. Since we will work only with the diagrams that are
closure of a braid, we will denote the singular resolution (one with the thick edge) by $\bar{E}_i$ if it comes from the
crossing of the type $\sigma_i$, and the obtained total resolution of $D$ will be the closure of the sequence of
$\bar{E}_i$'s.\\

{{
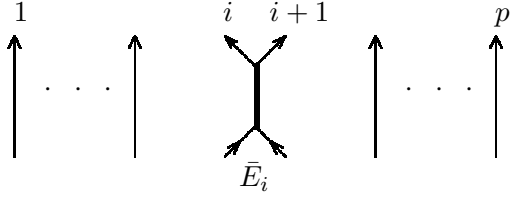
\begin{center}
\setlength{\unitlength}{4mm}
\begin{picture}(30,6) 
\linethickness{0.6pt}
%\line(1.00,2.00,1.00,3.00)
\qbezier(5.00,1.00)(5.00,3.00)(5.00,5.00)
\qbezier(4.80,4.60)(4.90,4.80)(5.00,5.00)
\qbezier(5.20,4.60)(5.10,4.80)(5.00,5.00)

\qbezier(9.00,1.00)(9.00,3.00)(9.00,5.00)
\qbezier(8.80,4.60)(8.90,4.80)(9.00,5.00)
\qbezier(9.20,4.60)(9.10,4.80)(9.00,5.00)

\linethickness{2.0pt}
\qbezier(13.00,2.00)(13.00,3.00)(13.00,4.00)
\linethickness{0.6pt}
\qbezier(12.00,1.00)(12.50,1.50)(13.00,2.00)
\qbezier(14.00,1.00)(13.50,1.50)(13.00,2.00)
\qbezier(12.00,5.00)(12.50,4.50)(13.00,4.00)
\qbezier(14.00,5.00)(13.50,4.50)(13.00,4.00)

\qbezier(14.00,5.00)(13.80,4.90)(13.60,4.80)
\qbezier(14.00,5.00)(13.90,4.80)(13.80,4.60)

\qbezier(12.00,5.00)(12.20,4.90)(12.40,4.80)
\qbezier(12.00,5.00)(12.10,4.80)(12.20,4.60)

\qbezier(12.50,1.50)(12.30,1.40)(12.10,1.30)
\qbezier(12.50,1.50)(12.40,1.30)(12.30,1.10)

\qbezier(13.50,1.50)(13.70,1.40)(13.90,1.30)
\qbezier(13.50,1.50)(13.60,1.30)(13.70,1.10)

\qbezier(17.00,1.00)(17.00,3.00)(17.00,5.00)
\qbezier(16.80,4.60)(16.90,4.80)(17.00,5.00)
\qbezier(17.20,4.60)(17.10,4.80)(17.00,5.00)

\qbezier(21.00,1.00)(21.00,3.00)(21.00,5.00)
\qbezier(20.80,4.60)(20.90,4.80)(21.00,5.00)
\qbezier(21.20,4.60)(21.10,4.80)(21.00,5.00)

\put(6.00,3.00){$\cdot$}
\put(7.00,3.00){$\cdot$}
\put(8.00,3.00){$\cdot$}
\put(18.00,3.00){$\cdot$}
\put(19.00,3.00){$\cdot$}
\put(20.00,3.00){$\cdot$}

\put(5.00,5.50){1}
\put(11.95,5.50){$i$}
\put(13.50,5.50){$i+1$}
\put(21.00,5.50){$p$}
\put(12.50,0.00){$\bar{E}_i$}

\end{picture}
\end{center}
}}

\indent Now, to each such resolution $D_{\epsilon}$ is assigned a graded vector space $\bar{V}_{\epsilon}\{-|\epsilon|\}$ (the theory is defined over $\Q$), for details see \cite{kovroz}, and then after summing along the columns of the cubic complex we obtain the chain groups ${C}_n^i(D)$ by: 
$${C}_n^i(D)=\bigoplus_{|\epsilon|=i}{\bar{V}_{\epsilon}}.$$
 Further, the differential $d_n$ is obtained as a signed sum of ``per-edge" maps, and thus we obtain the chain complex $(C_n(D),d_n)$. Denote its homology groups by $H_n^{i,j}(D)$. Finally, after overall shift we obtain the chain complex $\bar{C}_n(D)$ given by $$\bar{C}_n(D)=C_n(D)[-n_-(D)]\{-(n-1)\cdot n_+(D)+n\cdot n_-(D)\}.$$ 
\begin{remark}
There is also slight difference in $q$-gradings of the $sl(n)$ and standard $sl(2)$ theory. Namely, if we put  $n=2$ in $sl(n)$ theory (\cite{kovroz}) we obtain standard Khovanov homology (\cite{kov}) with $q$-grading inverted, i.e. the $q$-gradings of \cite{kovroz} are the negative of the $q$-gradings of \cite{kov}. We would obtain the same convention for the standard Khovanov homology (Section \ref{intro}) if we define $\deg 1=-1$, $\deg X=1$ and all later shifts in $q$-gradings $\{i\}$  we replace by $\{-i\}$. These conventions for $sl(2)$ theory are used in \cite{tan} and \cite{pat}.
\end{remark}

\indent Denote the homology groups of the complex $\bar{C}_n(D)$ by $\H^{i,j}_n(D)$. Then we have:
\begin{theorem}(\cite{kovroz})\label{KR}
The homology groups $\H^{i,j}_n(D)$ are independent of the planar projection $D$. Even more, the graded Euler characteristic of the complex $\C_n^{i,j}(D)$ is equal to $n$-specialization of HOMFLYPT polynomial.
\end{theorem}
Hence we can write $\H_n^{i,j}(K)$.\\

\begin{theorem}
For every positive braid knot $K$, we have that the homology group $\H^1_n(K)$ is trivial.
\end{theorem}
\textbf{Proof:}\\
\indent Mainly we will adapt our proof from the previous Section for this case. First of all, we again take the diagram $D$ of knot $K$, which is the closure of positive braid $\sigma$. Since $D$ has only positive crossings we have that $\H^1(K)$ is trivial if and only if $H^1(D)$ is trivial. So, we have to prove that the latter group is trivial.\\
\indent Again, we use the direct sum definition of  the chain groups $\bar{C}_n^0(D)$, $\bar{C}_n^1(D)$, $\bar{C}_n^2(D)$, and of the differentials $d_n^0$ and $d_n^1$.  We will prove that for each $t' \in \bar{C}_n^1(D)$ such that $d_n^1(t')=0$, there exists $t \in \bar{C}_n^0(D)$ such that $d_n^0(t)=t'$.\\
\indent We denote the restriction of $t'$ to the space $\bar{V}_{i\alpha}$ by $t'_{i\alpha}$. Again the differential ${d}_n^1$ maps only the spaces $\bar{V}_{i\alpha}$ and $\bar{V}_{i\beta}$ to the $\bar{V}_{i\alpha, i\beta}$, and  since $d_n^1(t')=0$ we have that its restriction to the $\bar{V}_{i\alpha, i\beta}$ is zero.\\
\indent From the definition we have that $\bar{V}_{i\alpha}=\bar{V}(\bar{E}_i)\{-1\}$ and $\bar{V}_{i\alpha, i\beta}=\bar{V}(\bar{E}_i^2)\{-2\}.$ From \cite{kovroz} we have that the last vector space is isomorphic to $\bar{V}(\bar{E}_i)\{-1\}\oplus\bar{V}(\bar{E}_i)\{-3\}$, and the projection of the per-edge map from $\bar{V}_{i\alpha}$ to $\bar{V}_{i\alpha, i\beta}$ onto the first summand of the latter space is identity map. \\
\indent Hence, we have that $t'_{i\alpha}=t'_{i\beta}$ for every $i$, $\alpha$, $\beta$, and our problem reduces, like in $sl(2)$ case, to the problem when there is at most one crossing of the type $\sigma_i$ for every $i$. In other words we are left with proving that $H^1_n(D')$ is trivial. However, since $D'$ is isotopic to the unknot the triviality of 
$H^1_n(D')$ follows from the independence of the $sl(n)$ homology of knot projection chosen, which concludes our proof. \kraj\\


\footnotesize
\begin{thebibliography}{99}

\bibitem{bn}
 D. Bar-Natan: \textit{On Khovanov's Categorification of the Jones Polynomial},  Alg. Geom. Top. 2: 337-370 (2002)
\bibitem{bn1}
 D. Bar-Natan: \textit{Khovanov's Homology for Tangles and Cobordisms}, arXiv:math.GT/0410495
\bibitem{bnprog}
 D. Bar-Natan: \textit{The Knot Atlas}, www.math.toronto.edu/~drorbn/KAtlas
\bibitem{dgr}
 N. Dunfield, S. Gukov and J. Rasmussen: \textit{The Superpolynomial for knot homologies}, arXiv:math.GT/0505662.
\bibitem{homfly}
 P.Freyd, D.Yetter, J.Hoste, W.B.R.Lickorish, K.Millet and A.Ocneanu: \textit{A new polynomial invariant of knots and links}, Bull. AMS (N.S.) 12, no. 2, 239-246 (1985) 

\bibitem{xy}
R.Gompf, A.Stipsicz: \textit{4-manifolds and Kirby calculus}, Graduate Studies in Mathematics, vol.20, American Mathematical Society, Providence, RI,1999.

\bibitem{gr}
 L. Helme-Guizon and Y. Rong: \textit{A Categorification for the Chromatic Polynomial}, Alg. Geom. Top. 5: 1365-1388 (2005)

\bibitem{jones}
 V.F.R.Jones: \textit{A polynomial invariant for knots via von Neumann algebras}, Bull. AMS (N.S.) 12, no. 1, 103-111 (1985) 
\bibitem{kauf}
 L. Kauffman: \textit{Knots and Physics}, 3ed., World Scientific, 2001.
\bibitem{kov}
 M. Khovanov: \textit{A categorification of the Jones polynomial}, Duke Math. J. 101:359-426 (2000). 
\bibitem{kol}
 M Khovanov: \textit{Categorifications of the Colored Jones polynomial}, arXiv:math.QA/0302060.
\bibitem{kov3}
 M. Khovanov: \textit{sl(3) Link Homology}, Alg. Geom. Top. 4: 1045-1081 (2004), arXiv:math.QA/0304375 
\bibitem{pat}
 M. Khovanov: \textit{Patterns in knot cohomology I}, Experiment. Math. 12: (2003), no. 3, 365-374, arXiv:math.QA/0201306.
\bibitem{tan}
 M. Khovanov: \textit{A functor-valued invariant for tangles}, arXiv:math.QA/0103190. 
\bibitem{kovroz}
 M. Khovanov, L. Rozansky: \textit{Matrix Factorizations and link homology}, arXiv:math.QA/0401268. 
\bibitem{lee}
 E.S. Lee: \textit{The support of the Khovanov's invariants for alternating knots}, arXiv:math.GT/0201105. 
\bibitem{lee2}
 E.S. Lee: \textit{On Khovanov invariant for alternating links}, arXiv:math.GT/0210213. 
\bibitem{moy}
 H. Murakami, T. Ohtsuki and S. Yamada: \textit{HOMFLY polynomial via an invariant of colored plane graphs}, Enseign. Math. (2) 44 (1998), no. 3-4, 325-360.
\bibitem{ras}
 J. Rasmussen: \textit{Khovanov homology and slice genus}, arXiv:math.GT/0402131.
\bibitem{ras2}
 J. Rasmussen: \textit{Knot polynomials and knot homologies}, arXiv:math.GT/0504045.
\bibitem{pt}
 J. Przytycki and P. Traczyk: \textit{Invariants of links of Conway type}, Kobe J. Math. 4 (1988), no. 2, 115--139.
\bibitem{shum}
 A. Shumakovich: \textit{KhoHo: a program for computing Khovanov homology}, www.geometrie.ch/KhoHo/
\bibitem{mojtorus}
 M. Sto\v si\' c: \textit{Homological thickness of torus knot}, arXiv:math.QA/0511
\bibitem{mojnov}
 M. Sto\v si\' c: \textit{Stable $sl(n)$-link homology for torus knots}, in preparation
\bibitem{v}
 O. Viro: \textit{Remarks on the Definition of the Khovanov homology}, arXiv:math.GT/0202199.

\end{thebibliography}
\end{document}